\documentclass[12pt]{amsart}
\usepackage{amsmath}
\usepackage{amssymb}
\usepackage{latexsym}
\usepackage{amscd}
\usepackage{citesort}

\newdimen\AAdi%
\newbox\AAbo%
\def\AAk#1#2{\s_etbox\AAbo=\hbox{#2}\AAdi=\wd\AAbo\kern#1\AAdi{}}%
\def\AAr#1#2#3{\s_etbox\AAbo=\hbox{#2}\AAdi=\ht\AAbo\raise#1\AAdi\hbox{#3}}%
\font\tenmsb=msbm10 at 12pt \font\sevenmsb=msbm7 at 8pt
\font\fivemsb=msbm5 at 6pt
\newfam\msbfam
\textfont\msbfam=\tenmsb \scriptfont\msbfam=\sevenmsb
\scriptscriptfont\msbfam=\fivemsb

\newtheorem{thm}{Theorem}

\newtheorem{cor}{Corollary}
\newtheorem{rem}{Remark}

\newcommand{\ba}{\begin{array}}
\newcommand{\ea}{\end{array}}

\newcommand{\Section}[2]{\setcounter{equation}{0}
\allowdisplaybreaks
\section[#1]{#2}}

\begin{document}
\title
[SOME REMARKS ON A MINKOWSKI SPACE $(R^n, F)$] {SOME REMARKS ON A MINKOWSKI SPACE $(R^n, F)$}
\author
[T.Q. Binh]{Tran Quoc Binh}
\address[T.Q. Binh]
{Institute of mathematics \\ University of Debrecen\\ H-4010
Debrecen, P.O. Box 12\\ Hungary} \email{binh@science.unideb.hu}

\thanks{Key words:Finsler space; Minkowski space; Mean curvature; Totally umbilical hypersurface.}

\renewcommand{\subjclassname}{%
  \textup{2000} Mathematics Subject Classification}
\subjclass{Primary 53C60; Secondary 53C20 }\date{} \maketitle

\begin{abstract}
We consider a complete, totally umbilical hypersurface $M$ of Riemannian space $(\hat{R}^n, \hat{g})$ induced by a Minkowski space $(R^n, F)$. Under certain conditions we prove that $M$ is isometric to a "round" hypersphere of the $(n + 1)-$dimensional Euclidean space. We also prove that the Minkowski norm $F$ must be arised from an inner product if there exist a non-zero vector field, which is parallel according to Levi-Civita connection of the metric tensor $\hat{g}$.

\end{abstract}

\Section{Introduction }{Introduction }

Let $(N, F)$ be an n- dimensional Finsler space, then for any $x\in N$ the restriction $F_x$ of $F$ to the tangent space $T_xN$ is a Minkowski norm. That is,  $F_x : T_xN \rightarrow R$ is non-negative, smooth function on $T_xN - \{0\}$, positive homogenous of degree one and strongly convex. The last condition means that the functions $ g_{ij}(x, y) := \frac{1}{2}\frac{F^2(x, y)}{\partial y^i\partial y^j}$ give a positive definite quadratic form on $T_xN - \{0\}$, that is $ g := g_{ij}(x, y)dy^idy^j$ is a Riemannian metric on $T_xN - \{0\}$.

If for any $x\in N$, $g_{ij}(x, y)$ doesn't depend on $y$, then $g = g_{ij}(x, y)dy^idy^j$ is an inner product on $T_xN$, and $G := g_{ij}(x, y)dx^idx^j$ is a Riemannian metric on $N$. Of course, any Riemannian space is also Finslerian, but the inverse is not true. It is easy to see that  $G := g_{ij}(x, y)dx^idx^j$ is Riemannian metric on $N$ iff the Cartan torsion tensor $C_{ijk} := \frac{F}{2}\frac{\partial g_{ij}}{\partial y^k}$ vanishes.

When does the Minkowski norm $F_x$ arise from an inner product on $T_xN$? Or equivalently when does the Finslerian metric function $F(x, y)$ on $TN$ determine a Riemannian metric $G := g_{ij}(x, y)dx^idx^j$ ? This is a very intersting and important question. Using the  maximum principles, the famous Deicke's theorem says that this occurs if the mean Cartan torsion $A_k := g^{ij}C_{ijk}$ vanishes.
Since for any $x\in N$, the tangent space $T_xN$ is an $n-$dimensional vector space, so for simplicity in this paper we consider the $n-$dimensional Euclidean space and a Minkowski norm $F$ on it. So, let $(R^n,F)$ be a Minkowski space and $g_{ij}(y) := (\frac{1}{2}F^2(y) )_{y^iy^j}.$  Then $\hat{g} := g_{ij}(y)dy^i\otimes dy^j$ defines a Riemannian metric on $\hat{R}^n := R^n - \{0\}$ and $(\hat{R}^n, \hat{g})$ becomes a Riemannian space.  It is well known that when $F$ is induced by an inner product on $R^n, \ (\hat{R}^n, \hat{g})$ is flat (that is, the curvature tensor vanishes). Further more, (1) any proper totally umbilical hypersurface is a "round" hypersphere and (2) any constant vector field is parallel according to a canonical Levi-civita connection induced by $\hat{g}$.

 In the present  paper, concerning a proper totally umbilical hypersurface $M$  of a Riemannian space $(\hat{R}^n, \hat{g})$ we try to find the answer for the following two questions: 1. Under which conditions is $M$  isometric to a "round" sphere. 2. If there exists a non-zero constant, and parallel vector field according to $\hat{g}$, then under which conditions is the Minkowski norm $F$ arised from an Euclidean inner product? Our main results are:
\begin{thm}
 Let  $M$ be a proper totally umbilical hypersurface of a Riemannian space $(\hat{R}^n, \hat{g})$ induced by a Minkowski space $(R^n, F),$ and let $M$ be no  level set  $S(r) := \{y| F(y) = r\}$, $r> 0.$ Suppose that  the normal part of $R(X, Y)Z$, $X, Y, Z \in \mathfrak{X}(M)$ vanishes,i.e.  $(R(X, Y)Z)^{\bot} = 0.$ Then $M$ must be isometric to a "round" hypersphere of the $(n + 1)-$dimensional Euclidean space.
 \end{thm}

 \begin{thm}
   Let   $(\hat{R}^n, \hat{g})$ be  a Riemannian space induced by a Minkowski space $(R^n, F).$ If $n \geq 3,$ and $F $ is absolutely homogeneous, then $F$ is arised from an inner product if and only if  there exists a non-zero constant vector field $b$ which is parallel according to the Levi-civita connection  $\hat{\nabla}$ of $\hat{g}$, that is if $\hat{\nabla}_Xb = 0$.
   \end{thm}

 \Section{Preliminaries}{Preliminaries}

In this section, we collect some facts  of the Minkowski space $(R^n, F)$. For more details  see
\cite{BCS}.  We use the Einstein convention. That is, repeated indexes with one upper  and one
lower indexes denote summation over their range  throughout this paper.

Let $(R^n,F)$ be a Minkowski space, and let $ (y^1, y^2, ..., y^n)$ denote a canonical global coordinate system on $R^n$. Then at every pont $y\in T_yR^n$ we have $\{\frac{\partial}{\partial y^i}\}$ and $\{dy^i\}$ which are the bases of the tangent space $T_yR^n$ and the cotangent space $T^{*}_yR^n$, respectively. The Minkowski norm on $R^n$ is a continuous function $F: R^n \rightarrow R^+$ which has the ptoperties:\\
\begin{enumerate}
  \item  $F$ is smooth on  $\hat{R}^n := R^n - \{0\}$
  \item $F(y) > 0$  for all  vector $y\in \hat{R}^n$
  \item positive homogenous, i.e., $F(\lambda y) = \lambda F(y)$ for all $y\in R^n$ and all $\lambda > 0$, and
    \item  $F$is strongly convex, i.e., the quantities $g_{ij}(y) := \frac{1}{2}\frac{F^2}{\partial y^i\partial y^j}$ form a positive definite matrix. Thus in this case $\hat{g}(y) := g_{ij}(y)dy^idy^j$ defines a metric tensor for $ \hat{R}^n$
\end{enumerate}

Now we recall some important facts related to the Riemannian space $(\hat{R}^n, \hat{g})$ (\cite{BCS}): \\
Let $\hat{\nabla}$ and $\hat{\gamma}_{jk}^i$ denote the Levi-Civita connection and the Christoffel symbols of the second kind of the metric tensor  $\hat{g}.$ Then
$$\hat{\nabla}_{\frac{\partial}{\partial y^k}}\frac{\partial}{\partial y^j}:= \hat{\gamma}_{jk}^i\frac{\partial}{\partial y^i},\eqno{(2.1)}$$
where
$$ \hat{\gamma}_{jk}^i = \frac{g^{is}}{2}(\frac{\partial g_{sj}}{\partial y^k} + \frac{\partial g_{sk}}{\partial y^j} - \frac{\partial g_{jk}}{\partial y^s}) . \eqno (2.2)$$
It is easy to see that
$$ \hat{\gamma}_{jk}^i = \frac{1}{F}C_{jk}^i, \eqno (2.3)$$
where
$$C_{jk}^i(x,y):=g^{ks}C_{ijs}(x,y)$$ comes from the Cartan torsion
$$C_{ijk}(x,y):=\frac{F}{4}\frac{\partial^3F^2(x,y)}{\partial
y^i\partial y^j\partial y^k}$$
The curvature tensor of $\hat{\nabla}$ is
$$\hat{R}(U, V)W : = (\hat{\nabla}_U\hat{\nabla}_V - \hat{\nabla}_V\hat{\nabla}_U - \hat{\nabla}_{[U, V]})W, \ U, V, W \in \mathfrak{X}(M).$$
In a local coordinate system,
$$ \hat{R}(U, V)W = W^jU^kV^l\hat{R}_{j jkl}^i\frac{\partial}{\partial y^i}, $$
where  in our case, it can be shown that
$$\hat{R}_{j kl}^i = \frac{1}{F^2}(C_{j k}^sC_{s l}^i - C_{j l}^sC_{s k}^i).$$

For any two linearly independent tangent vectors $U, V$ in $T_y\hat{R}^n$ the sectional curvature corresponding to the 2-plane defined by $U$ and $V$ is given by
$$\hat{K}(U, V) := \frac{\hat{g}(\hat{R}(V, U)V, V)}{\hat{g}(U, U)\hat{g}(V, V) - (\hat{g}(U, V))^2}.$$
If the sectional curvature $K$ does not depend on the choice of the 2-plane, then we speak  of constant curvature or the space is said to be a space form.\\

In the case of a level set $S(r) := \{y\in R^n: F(y) = r> 0\},$ we have
\begin{thm}{(\cite{BCS},Prop. 14.6.1, p. 401)}
The sectional curvature $K(U, V)$ of a level  hypersurface $S(r)$ is related to the sectional curvature $\hat{K}(U, V)$ of $(\hat{R}^n, \hat{g})$ by \\
 $K(U, V) = \hat{K}(U, V) + \frac{1}{r^2}$ \\
  and the following statemants are equivalent:\\
(a) The Riemannian space $(\hat{R}^n, \hat{g})$ is flat.\\
(b) For any $r > 0$, the level  hypersurface $S(r)$ has constant sectional curvature $\frac{1}{r^2}$\\
(c) For some $r_0 > 0$, the level  hypersurface $S(r_0)$ has constant sectional curvature $\frac{1}{r_0^2}$\\
\end{thm}
In the proof of our results we will also need the following important theorems.
\begin{thm}{( \cite{BCS}, Theorem 14.9.2 (Brickell-Theorem), p.415)} Let $(R^n, F)$ be a Minkowski space. If $n \geq 3$, $F$ is absolutely homogeneous of degree 1 and the Riemannian space $(\hat{R}^n, \hat{g})$ is flat, then $F$ must be the norm of an inner product on $R^n$.
\end{thm}
From Theorem 3 and Theorem 4 we have
\begin{cor}{(\cite{BCS})}
If $n \geq 3$, $F$ is absolutely homogeneous of degree 1 and for some $r > 0$, $(S(r), \hat{g})$ has constant sectional curvature $K = \frac{1}{r^2}$, then the Minkowski norm is induced by an inner product.
\end{cor}
\begin{thm}{(Obata's Theorem, \cite{O})}
In order for a complete Riemannian manifold of dimesion $n \geq 2$ to admit a non-constant function $\varphi$ with $\nabla_Xd\varphi + c^2 \varphi X = 0$ for any vector filed $X$, it is necessary and sufficient that the manifold be isometric with a sphere $S^n(c)$ of radius $\frac{1}{c}$ in the $(n + 1)$-Euclidean space.
\end{thm}
Obata's theorem is very useful to get important results in  geometry of submanifolds, firstly based on  the works of Bang-Yen Chen (\cite{Ba}) and its version for Finsler geometry is given by B. Badabad (\cite{Bi}). We can find the studies of submanifolds of Minkowski spaces in work of X. Cheng and J. Yan (\cite{Ch-Ya}), and the studies of umbilical hypersurfaces of Minkowski spaces in the sense of Finsler geometry can be found in the paper of
J. Li (\cite{Li}).
\Section{Proof of Theorem 1.}{Proof of Theorem 1.}
For the totally umbilical hypesurface $M$ of Riemannian space $(\hat{R}^n, \hat{g})$  the Gauss equation is
$$\hat{\nabla}_XY = \nabla_XY + h(X, Y),$$
where $h$ is the second fundamental form.
In our case $(M, g)$ is a proper totally umbilical hypersurface. So we have
$$\hat{\nabla}_XY = \nabla_XY + Hg(X, Y)\nu ,$$
where $H$ and $\nu$ denote the mean curvature  and the unit normal vector field of the hypersurface $M$ and $H \neq 0$. We choose the unit normal vector $\nu$ in such a way that the mean curvatute  $H$ takes positive value somewhere on $M$. From the identity $g_{ij}(y)y^iy^j = F^2(y)$ it follows that $\frac{y}{F(y)}$ is the unit normal vector field of the level hypersurface $S(r)$.   Since $F^2(y)$ is positive homogenous of degree 2, the functions $g_{ij}(y)$ are positive homogenous of degree 1, from $(2.2)$ we have
 $$\hat{\gamma}_{ij}^k(y)y^j = C_{ij}^k(y)y^j = 0,$$
 from which we obtain
$$\hat{\nabla}_Xy = X, \ \forall X\in \mathfrak{X}(\hat{R}^n).$$
Let us define a function $f(y) := \hat{g}(y,\nu)$ on $M$. Since $M$ is not a level set, the function $f$ is non-constant. For, if $f(y)$ is constant, then for any vector field $X$ on $M$, we have\\
$0 = X(\hat{g}(y, \nu)) = \hat{g}(\hat{\nabla}_Xy, \nu) + \hat{g}(\hat{\nabla}_X\nu, y) = \hat{g}(X, \nu) + \hat{g}(y, - HX) = -H\hat{g}(y, X)$,\\
from which $X(F^2(y)) = X(\hat{g}(y, y)) = 2\hat{g}(\nabla_Xy, y) = \hat{g}(X, y) = 0$. This means $F(y) = constans.$ Its contradiction to the fact that $M$ is not a level set. \\
 Now we compute the gradient, $\nabla f$ of a function $f$ in $M$. We have
$$X(f) = X\hat{g}(y, \nu) = \hat{g}(\hat{\nabla}_Xy, \nu) + g(y, \hat{\nabla}_Y\nu)$$
$$= \hat{g}(X, \nu) + \hat{g}(y, - HX) = - H\hat{g}(y, X) = - H\hat{g}(y^{T}, X) = - Hg(y^T, X),$$
where $y^T := y - \hat{g}(y, \nu)\nu$ is the tangential part  to $M$ of $y.$ From $X(f) = g(X, \nabla f),$ it follows
$$\nabla f = - H y^{\top}.$$
Now, we prove that under our assumptions the mean curvature $H$ is constant. Recall the Ricci equation from the geometry of submanifold
$$(\hat{R}(X, Y)Z)^{\bot} = (\nabla h)(X, Y, Z)  - (\nabla h)(Y, X, Z),$$
where
$$(\nabla h)(Y, X, Z) := \nabla_X^{\bot}h(Y, Z) - h(\nabla_XY, Z) - h(Y, \nabla_XZ).$$
At any $y\in M$ we can choose a local normal coordinate system such that $\nabla _XY(y) = \nabla_YX(y) = 0.$ Then we obtain
$$\nabla_Y^{\bot}H\hat{g}(X, Z)\nu = \nabla_X^{\bot}H\hat{g}(Y, Z)\nu .$$
Choose the unit vector fields $X, Y, Z$ such that $X = Z$ and orthogonal to $Y.$ Then we have $\nabla_Y^{\bot}H\nu =0,$ from which it follows that $Y(H) = 0$, for any tangent vector $Y \in T_yM$. This means the mean curvature $H$ is constant.\\
Now,
$$\nabla_X\nabla f = (\hat{\nabla}_X\nabla f)^{\top} = - H(\hat{\nabla}_Xy^{\top})^{\top} =$$
$$= -H\hat{\nabla}_X(y - \hat{g}(y, \nu)\nu) = -H\{X + \hat{g}(y, \nu)HX\} = -H^2(\frac{1}{H} + f(y))X.$$
Let now $\tilde{f}(y) := f(y) + \frac{1}{H}.$ Then we have
$$\nabla_X\nabla \tilde{f} + H^2\tilde{f}X = 0.$$ From the theorem of Obata it follows that the totally umbilical hypersurface $M$ is isometric to a round sphere $S^n(c)$ of radius $\frac{1}{H}$ in the $(n + 1)$-Euclidean space.$\Box$
\begin{rem}
We know that the level hypersphere $S^n(r)$ is totally umbilical with unit normal vector field $\frac{y}{r}$ at any $y\in S^n(r). $  It is easy to show that the converse is also true. For, if $\nu = \frac{y}{F(y)}$ is a unit normal vector field at a point $y$ of a totally umbilical hypersurface $M$, then for any tangent vector field $X$, we have
 $$ 0 = X\hat{g}(\nu, \nu) = X \hat{g}(\frac{y}{F(y)}, \frac{y}{F(y)}) = 2\hat{g}(X(\frac{1}{F})y + \frac{y}{F(y)}X, y) = 2X(\frac{1}{F})F^2(y).$$
Since $M$ is a submanifold of $\hat{R}^n$, $F(y) > 0$ for any $y \in M$, thus $X(\frac{1}{F}) = 0$, and we obtain $F(y) = r =$ constant, $M$ is nothing but the level set $S(r).$
\end{rem}

\section{Proof of Theorem 2.}
We suppose that $b$ is a non-constant vector filed, which is parallel according to the Levi-Civita connection $\hat{\nabla}$ of a Riemannian space $(\hat{R}^n, \hat{g})$ i. e. $\hat{\nabla}b = 0.$ We consider the level set $S(1) := \{y\in \hat{R}^n: F(y) = 1\} ,$ and we show that $S(1)$ is a hypersurface of $(\hat{R}^n, \hat{g})$ whose sectional curvature is constant 1. Then, by Theorem 3. the Riemannian space $(\hat{R}^n, \hat{g})$ is flat and then according to Theorem 4. $F$ must be a norm induced by an inner product on $R^n.$\\
Let us define  $f(y) := \hat{g}(y, b)$. For any tangent vector field $X$ on $S(1)$ we have
$$X(f) = X\hat{g}(y, b) = \hat{g}(\hat{\nabla}_Xy, b) - \hat{g}(\hat{\nabla}_Xb) =$$
$$= \hat{g}(X, b) = \hat{g}(X, b^{\top)}.$$
Thus, we obtain the gardient of the function $f$ on $S(1)$,
$$\nabla f = b^{\top} = b - \hat{g}(b, y)y.$$
Now,
$$\nabla_X\nabla f = (\hat{\nabla}_Y\nabla f)^{\top} = (\hat{\nabla}_Y(b - \hat{g}(b, y)y))^{\top} =$$
$$= \{ - X\hat{g}(b, y)y - \hat{g}(b, y)\hat{\nabla}_Yy\}^{\top} = -f(y)X.$$
Since  $b$ is constant vector field its easy to see that $f(y)$ is non-constant function on $S(1).$ From the Obata's theorem (Theorem 5), it follows that $S(1)$ is isometric to the round sphere of radius 1, and thus its sectional curvature is also constant 1. This completes the  proof of Theorem 2.
$\Box$

\end{document}